\documentclass[12pt,letterpaper]{article}
 \usepackage[T1]{fontenc}
 \usepackage{amsfonts}
 \usepackage{amssymb}
 \usepackage{stmaryrd}
 \usepackage{amsthm,amsmath}
 \usepackage{a4wide}
 \usepackage{ulem}
 \usepackage{enumerate}
 \usepackage[dvips]{graphicx}
 \usepackage{color}
 \def\sym{\fam\comfam\com}
 \font\tensym=msbm10
 \font\sevensym=msbm7
 \font\fivesym=msbm5
 \newfam\symfam
 \textfont\symfam=\tensym
 \scriptfont\symfam=\sevensym
 \scriptfont\symfam=\fivesym
 \def\sym{\fam\symfam\relax}
 
 \def\R{{\sym R}}

\begin{document}
\bibliographystyle{plain}
 \title{\bf{On the existence of periodic solution of perturbed 
generalized
Li\'enard equations }} \date{}
 \maketitle
 \newtheorem{intro}{introduction}
 \newtheorem{theo}{Theorem}
 \newtheorem*{nt}{Remark}
 \newtheorem{prop}{Proposition}
 \newtheorem{cor}{Corollary}
 \newtheorem*{pro}{proof}
 \renewcommand\abstractname{Abstract}
 \centerline{by}
\par\medbreak
\centerline{Islam Boussaada}
\centerline{LMRS, UMR 6085, Universite de Rouen}
\centerline{Avenue de l'universit\'e, BP.12}
\centerline{76801 Saint Etienne du Rouvray, France.}
\centerline{email:islam.boussaada@etu.univ-rouen.fr}
\par\medbreak
\centerline{and}
\par\medbreak
\centerline{A. Raouf Chouikha}
\centerline{Universite Paris 13 LAGA}
\centerline{Villetaneuse 93430, France.}
\centerline{email: chouikha@math.univ-paris13.fr}
\par\medbreak
\begin{abstract}

Under conditions of Levinson-Smith we prove the existence 
of a $\tau$-periodic solution for the perturbed generalized Li\'enard equation
 $$u''+\varphi(u,u')u'+\psi(u)=\epsilon\omega(\frac{t}{\tau},u,u')$$
with periodic forcing term. \\
We deduce sufficient condition for existence of a periodic solution for the equation
$$ u''+\sum_{k=0}^{2s+1} p_k(u){u'}^k=\epsilon\omega(\frac{t}{\tau},u,u').$$
In particular, our method can be applied to the equation
$$u''+[u^2+(u+u')^2-1]u'+u=\epsilon\omega(\frac{t}{\tau},u,u').$$
Finally, these results will be illustrated by some numerical exemples.
\begin{center} {\it{Keywords : perturbed systems, Li\'enard equation, periodic solution }} \end{center}

\end{abstract}
\section{Introduction }
 Consider Li\'enard equation \begin{equation}
  u''+\varphi(u)u'+\psi(u)=0
                             \end{equation}
where $u'=\frac{du}{dt}$ , $u''=\frac{d^2u}{dt^2}$, $\varphi$  and  
$\psi$ are $C^{1}. $
 Existence problem of periodic solution of period $\tau_0$ was the purpose's study of 
many
authors. M.Farkas selected some typical works on this subject see [3],
where the Poincar\'e-Bendixson theory plays a crucial role.\\
In general, a periodic perturbation of the Li\'enard equation does 
not possess a periodic solution as described the Moser example, [1].\\

Let us consider the perturbed Li\'enard equation of the form 
\begin{equation}
u''+\varphi(u)u'+\psi(u)=\epsilon\omega(\frac{t}{\tau},u,u')
\end{equation}
where $\omega$ is a $controllably\; periodic\; perturbation$ in the 
Farkas
sense, i.e. it is periodic with a period $\tau$ which can be choosen 
appropriately. The existence of a non trivial periodic solution for (2) was
studied by Chouikha [1]. Under very mild conditions it is proved that
to each small enough amplitude of the perturbation there belongs a one
 parameter family of periods $\tau$ such that the perturbed system has a 
unique
 periodic solution with this period. \\
 
 Let us consider now the following 
generalized
 Li\'enard equation  that is "a more realistic assumption 
in modelling many real world phenomena" ( [3] page 105)
 \begin{equation}\label{l}
 u''+\varphi(u,u')u'+\psi(u)=0.
 \end{equation}
Where
 $\varphi$  and  $\psi$ are $C^{1} $ and satisfy some assumptions 
that will
 be specified below.\\
 The leading work of investigation for the existence of periodic
solution of generalized Li\'enard systems
was established by Levinson-Smith [4].\\
 Let us define conditions $C_{LS}$ .\\
 
 {\bf Definition}\quad {\it The functions $\varphi$ and  $\psi$ satisfy $C_{LS}$ if :
 $$x\psi(x)>0\qquad\ for \qquad|x|>0.$$
 $$\int_0^{x} \psi(s) ds =\Psi(x) \quad and \quad 
lim_{x\rightarrow +\infty}
 \Psi(x) = + \infty , \quad \varphi(0,0) < 0.$$

Moreover, there exist some numbers $0 < x_0 < x_1$ and $M>0$ such that :

\begin{eqnarray*}
 \varphi(x,y)&\geq& 0 \qquad\;\;\; \;\;for\qquad |x|\geq 
x_0\\
\varphi(x,y)&\geq& -M \qquad for \qquad |x|\leq  x_0\\
x_1 > x_0 ,\qquad \int_{x_0}^{x_1} \varphi(x,y(x))dx\quad &\geq& 
\quad 10 M x_0
\end{eqnarray*}
for every decreasing function $y(x)>0.$}

\bigskip

\begin{prop}\text{ (Levinson-Smith) }\\
{\it{
When the functions $\varphi$ and  $\psi$ are of class $C^1$ and satisfy $C_{LS}$ then the generalized Lienard equation (\ref{l}) has at least one non-constant 
$\tau_0$-periodic solution.}}
\end{prop}
A non trivial solution will be denoted $u_0(t)$, and its period $\tau_0$.
This proposition has many improvements (under weaker hypotheses) due to Zheng-Zuo-Huan,
 N.Wax-P.J. Ponzo among other authors, [3]. \\
 
 Our paper is organized as follows.\\
At first, we prove the existence of a periodic solution for the perturbed generalized Li\'enard equation
\begin{equation}\label{l1}
u''+\varphi(u,u')u'+\psi(u)=\epsilon\omega(\frac{t}{\tau},u,u'),
\end{equation}
Where  $t,\epsilon,\tau\in \sym R$ are such that 
$|\tau-\tau_0|<\tau_1<\tau_0$, $|\epsilon|<\epsilon_0$ with $\epsilon_0\in \sym R$ 
sufficient small and $\tau_1$ is a fixed real scalar.
We will use the Farkas method which 
itself was effective for perturbed Li\'enard equation.\\
In the third section, we will propose a criteria for the existence
 of periodic solution for 
\begin{equation}\label{l2}
u''+\sum_{k=0}^{2s+1} p_k(u){u'}^k=\epsilon\omega(\frac{t}{\tau},u,u'),
\end{equation}
with $ s\in \sym N $ and $p_k$ are $C^1$ functions, $\forall k\leq2s+1.$\\
In the second part of the section,  using a result of De Castro ([2]) 
 we will prove uniqueness of a periodic solution for the equation
\begin{equation}\label{l4}
u''+[u^2+(u+u')^2-1]u'+u=0.
\end{equation}
Sufficient condition of the existence of periodic solution for the equation
\begin{equation}\label{l3}
u''+[u^2+(u+u')^2-1]u'+u=\epsilon\omega(\frac{t}{\tau},u,u').
\end{equation}
will be found.
At the end of the paper, some phase plane exemples are given in order to illustrate the above results. In particular, we describe uniqueness of a solution for equation (6) and the existence of a solution of equation (7) for \ $\omega(\frac{t}{\tau},u,u') = (\sin 2t) \ u'$.

\section{Periodic solution of perturbed generalized Lienard equation }
In this part of this paper we will deal with the proof of existence of 
periodic solution of the pertubed generalized Lienard equation(\ref{l1}) such that
 the unperturbed one (\ref{l}) has at least one periodic solution.\\

The method of proof that we will employ was described in [1] and [3].\\

Consider the equation (\ref{l})
$$u''+\varphi(u,u')u'+\psi(u)=0.$$

We assume that $\varphi$ and $\psi$ are $C^1$ and satisfy $C_{LS},$ 
then by Proposition 1 there exists at least a non trivial periodic solution denoted $u_0(t).$ 

Let the least positive period of the solution $u_0(t)$ be denoted 
by $\tau_0$ and $U \textrm { an open subset of } \R^2$ containing $(0,0).$\\
These notations will be used in the sequel of the paper.

\begin{theo}\label{monpre}
Let $\varphi$ and $\psi$ be $C^1$ and satisfy $C_{LS}.$\\
Suppose $1$ is a simple characteristic multiplier of the variational system 
associated to (\ref{l}).
 Then there are two real functions
$\tau,h$ defined on $U\subset\R^2$ and constants $\tau_1<\tau_0 $
 such that the periodic solution 
$\nu(t,\alpha,a+h(\epsilon,\alpha),\epsilon,\tau(\epsilon,\alpha))$
  of the equation
 $$u''+\varphi(u,u')u'+\psi(u)=\epsilon\omega(\frac{t}{\tau},u,u'),$$
 exists for $(\epsilon,\alpha)\in U $, $|\tau-\tau_0|<\tau_1$, 
 $\tau(0,0)=\tau_0$ and $h(0,0)=0.$
\end{theo}
We point out that the characteristic multipliers are the eigenvalues of 
the characteristic
 matrix which is the fundamental matrix in the time $\tau_0.$\vspace{0.3cm}\\
{\textbf{Proof of Theorem \ref{monpre}:}}\\ 
Following the method used in [3], 
we admit this setting 
$x_2=u$ , $x_1=\frac{du}{dt}=u'$ \\
and note $x=col(x_1,x_2)= col(u',u).$\\
The plane equivalent system of (\ref{l}) is:
\begin{equation}\label{l5}
x'=f(x)\Longleftrightarrow \left\{\begin{array} 
{rl}x'_1&=-\varphi(x_2,x_1)x_1-\psi(x_2)\\x'_2&=x_1
                                  \end{array}\right. \end{equation}

with $$f(x)=col(-\varphi(x_2,x_1)x_1-\psi(x_2),x_1).$$
Then the system  (\ref{l5}) has the periodic solution $q(t)$ with period 
$\tau_0.$\\
We define  
$$q(t)=col({u_0}'(t),u_0(t))$$ and therefore
 $$q'(t)=col(-\varphi(u_0(t),{u_0}'(t)){u_0}'(t)-\psi(u_0(t)),{u_0}'(t)) .$$
The variational system  associated to (\ref{l5}) is  
\begin{equation}\label{l6}
y'={f_x}'(q(t))y,
\end{equation}

 Without loss of generality, we take the initial conditions
 $$ t=0 ,\qquad u_0(0)=a <0 \,\,and\,\,{u_0}'(0)=0$$ Hence
  $${f_x}'(q(t))=\left
  (\begin{array}{cc}-{\varphi'}_{x_1}(u_0(t),{u_0}'(t)){u_0}'(t)-\varphi(u_0(t),{u_0}'(t))
&-{\varphi'}_{x_2}(u_0(t),{u_0}'(t)){u_0}'(t)-\psi'(u_0(t))\\1&0\end{array}\right)$$

Notice that 
$q'(t)=col(-\varphi(u_0(t),{u_0}'(t)){u_0}'(t)-\psi(u_0(t),{u_0}'(t))$
 is the first solution of the variational system. Now
 we calculate the second one, denoted 
$\widehat{y}(t)=col(\widehat{y}_1(t),\widehat{y}_2(t))$
 linearly independent with $q'(t)=y(t),$ in order to write the 
fundamental matrix.\\
 Consider 
$$I(s)=exp[-\int_0^s({\varphi'}_{x_1}(u_0(\rho),{u_0}'(\rho)){u_0}'(\rho)+
 \varphi(u_0(\rho),{u_0}'(\rho)))d\rho]$$ and
$$\pi(t)=-\int_0^t(\varphi(u_0(\rho),{u_0}'(\rho)){u_0}'(\rho)+\psi(u_0(\rho))^{-2}
 ({\varphi'}_{x_2}(u_0(t),{u_0}'(t)){u_0}'(t)+\psi'(u_0(t)))I(\rho)d\rho$$
 We then obtain 
$$\widehat{y}_1(t)=-[\varphi(u_0(t),{u_0}'(t)){u_0}'(t)+\psi(u_0(t)]\pi(t)$$
 $$\widehat{y}_2(t)={u_0}'(t)\pi(t)+\pi'(t)
 \frac{\varphi(u_0(t),{u_0}'(t)){u_0}'(t)+\psi(u_0(t)}
 {{\varphi'}_{x_2}(u_0(t),{u_0}'(t)){u_0}'(t)+\psi'(u_0(t))}$$
 It is known by [1] or [3] that the fundamental matrix satisfying 
$\Phi(0)=Id_2 $ is written as
  
$$\Phi(t)=\left(\begin{array}{cc}\frac{\varphi(u_0(t),{u_0}'(t)){u_0}'(t)+\psi(u_0(t))}{\psi(a)}
&\psi(a)\pi(t)[\varphi(u_0(t),{u_0}'(t)){u_0}'(t)+\psi(u_0(t)]
     \\-\frac{{u_0}'(t)}{\psi(a)}&-\psi(a)u'_0(t)\pi(t)-\psi(a)\pi'(t)
     \frac{\varphi(u_0(t),{u_0}'(t)){u_0}'(t)+\psi(u_0(t))}
{{\varphi'}_{x_2}(u_0(t),{u_0}'(t)){u_0}'(t)+\psi'(u_0(t))}\end{array}\right) $$

Thus, $$\Phi(\tau_0)= \left(\begin{array}{cc}1&{\psi(a)}^2 
\pi(\tau_0)\\
0&\rho_2\end{array}\right).$$
We use the Liouville's formula
 $$det\Phi(t)=det\Phi(0)exp\int_0^t Tr ({f_x}'(q(\tau)))d\tau.$$
  As $det(\Phi(0))=1, $
 we deduce the characteristic multipliers associted to (\ref{l6}):\\ 
$\rho_1=1$ and
 $\rho_2=I(\tau_0)=exp[-\int_0^{\tau_0}({\varphi'}_{x_1}(u_0(\rho),{u_0}'(\rho)){u_0}'(\rho)+
 \varphi(u_0(\rho),{u_0}'(\rho)))d\rho].$\\

From [3], we have : 
$$J(\tau_0)=-Id_2+\left[\begin{array}{cc}-\psi(a)&0
\\0&0\end{array}\right]+\Phi(\tau_0)$$
Hence we obtain the jacobian matrix :
$$J(\tau_0)=\left(\begin{array}{cc}-\psi(a)&{\psi(a)}^2 \pi(\tau_0)\\
0&\rho_2-1\end{array}\right),$$\\

since $1$ is a simple characteristic multiplier $(\rho_2\neq1),$ $$detJ(0,0,0,\tau_0)\neq0 .$$
We define the periodicity condition
\begin{equation}\label{lz} z(\alpha,h,\epsilon,\tau):=\nu(\alpha+\tau,a+h,\epsilon,\tau)-(a+h)=0
 \end{equation}
By the Implicit Function Theorem there are ${\epsilon}_0>0$ and ${\alpha}_0>0$ and uniquely 
determined functions $\tau$ and $h$ defined on
 $U=\{(\alpha,\epsilon)\in{\R}^{2}:|\epsilon|< \epsilon_0,|\alpha|<{\alpha}_0\}$ 
 such that :\,$\tau,h\in C^1,$ $\tau(0,0)=T_0,h(0,0)=0$ and $z(\alpha,h,\epsilon,\tau)\equiv0.$
 Since (\ref{lz}), the periodic solution of (\ref{l1}) is with the period $\tau(\epsilon,\alpha)$
  near $T_0$ and with path near the path of the unperturbed solution.
  \rightline{$\Box$} 
  \\
 
 In particular if $\rho_2<1$, the periodic solution is orbitally asymptotically stable i.e. stable 
 in the Liapunov sense and it is attractive see ([3] page 346).\\
 Thus, the following inequality is a criteria of the existence of orbital
  asymptotical stable periodic solution of the equation (\ref{l1}).
 
\begin{equation}\label{l7}
\rho_2<1\Longleftrightarrow
  \int_0^{\tau_0}({\varphi'}_{x_1}(u_0(\rho),{u_0}'(\rho)){u_0}'(\rho)+
 \varphi(u_0(\rho),{u_0}'(\rho)))d\rho>0. 
 \end{equation}
 
Using Proposition 1, we conclude the existence of 
non trivial periodic solution for perturbed generalized Li\'enard equation.
 \section{Results on the periodic solutions }
 \subsection{Special case }
Let us now consider the equation 
\begin{equation}\label{l8}
u''+\sum_{k=0}^{2s+1} 
p_k(u){u'}^k=0.
\end{equation}
Let $p_k$ be $C^1$ function, $\forall k\leq2s+1$ for $ s\in \sym N .$\\
This is a special case of Li\'enard equation with $$p_0(u)=\psi(u)$$  and  
$$\varphi(u,u')=\sum_{k=1}^{2s+1} p_k(u){u'}^{k-1} .$$
We will suppose $\varphi$ and $\psi$ verify $C_{LS}$ conditions. Let 
 $U $ be an open subset of $ \R^2$ containing $(0,0).$\\

The associated perturbed equation, as denoted previously,  is equation (\ref{l2}) 
$$ u''+\sum_{k=0}^{2s+1} 
p_k(u){u'}^k=\epsilon\omega(\frac{t}{\tau},u,u'),$$

\begin{nt}
The last non-zero term of the finite sum $\sum_{k=0}^{2s+1} p_k(u){u'}^k$ has an odd index.\\
Then it is necessary to have the element $x_0 \neq 0$ in the $C_{LS}$ conditions.\\
\end{nt}

\begin{theo}\label{monsec} Let $\varphi$ and $\psi$ be  $C^1$ and satisfy $C_{LS}$.\ 
If $1$ is a simple characteristic multiplier of the variational system 
associated to (\ref{l8})
 then there are two functions
$\tau,h :U\longrightarrow R$ and constants $\tau_1<\tau_0 $
 such that the periodic solution 
$\nu(t,\alpha,a+h(\epsilon,\alpha),\epsilon,\tau(\epsilon,\alpha))$
  of the equation
 $$ u''+\sum_{k=0}^{n} p_k(u){u'}^k=\epsilon\omega(\frac{t}{\tau},u,u')$$
 exists for $(\epsilon,\alpha)\in U $ with $|\tau-\tau_0|<\tau_1$, $\tau(0,0)=\tau_0$ and $h(0,0)=0.$
\end{theo}


{\textbf{Proof of Theorem \ref{monsec}:}}\\
We will use the same method and we will proceed as previously in the  
existence theorem
 of non-trivial periodic solution of the perturbed system.\\
Consider the unperturbed equation to compute some useful elements.
First we assume that $2s+1=n,$ to simplify the notations.
Let\;$x_2=u$ and $x_1=\frac{du}{dt}=u'.$\\
The equivalent plane system of (\ref{l8}) is

\begin{equation}\label{l9} 
x'=f(x)\Longleftrightarrow\left\{\begin{array} 
{rl}x'_1&=-\sum_{k=0}^{n} p_k(x_2){x_1}^k\\x'_2&=x_1
                                  \end{array}\right.
 \end{equation}

with $$f(x)=col(-\sum_{k=0}^{n} p_k(x_2){x_1}^k,x_1).$$
Let $q(t)=col(u'_0(t),u_0(t))$ the periodic solution of (\ref{l9}).

The variational system associated to (\ref{l9}) is $$y'=f'_x(q(t))y$$
with the periodic solution 
 $$q'(t)=col(-\sum_{k=0}^{n} p_k(u_0)(t){u_0'}^k(t) , u'_0(t)),$$
 
hence $$f'_x(q(t))=\left(\begin{array}{cc}-\sum_{k=1}^{n} k 
p_k(u_0(t)){u'_0(t)}^{k-1}
&-\sum_{k=0}^{n} p'_k(u_0(t)){u_0'(t)}^k\\1&0\end{array}\right) .$$
We assume the initial values :
$$ t=0 ,\qquad u_0(0)=a <0\,\, and\,\, {u_0}'(0)=0.$$
 Then $q(0)=col(0,a)$ and $q'(0)=col(-\psi(a),0).$\\

By the same way as the previous section we compute the fundamental 
matrix associated
to(\ref{l9}) denoted $ \Phi(t).$ 
Determine the second vector solution (linearly 
independent with $q'(t)=y(t)$).\\A trivial
calculation described in [1] and [3] gives us the second solution  
denoted $\widehat{y}(t),$
hence $\Phi(t)=(\frac{y(t)}{y(0)},y(0)\widehat{y}(t)).$
For that consider $$I(s)=exp[-\int_0^s(\sum_{k=1}^{n} k 
p_k(u_0(\rho)){u'_0(\rho)}^{k-1})d\rho],$$
and denote as in the previous section $$ \pi(t)=-\int_0^t(\sum_{k=0}^{n}
p_k(u_0)(\rho){u_0'}(\rho)^k)^{-2}(\sum_{k=0}^{n} 
p'_k(u(t)){u'}^k(t))I(\rho)d\rho.$$
Sine $\widehat{y}(t)=col(\widehat{y}_1(t),\widehat{y}_2(t)),$ where
$$\widehat{y}_1(t)=-(\sum_{k=0}^{n} p_k(u_0)(t){u_0'}(t)^k)\pi(t) $$
$$\widehat{y}_2(t)=u'_0(t)\pi(t)+\pi'(t)\frac{\sum_{k=0}^{n} 
p_k(u_0)(t){u_0'}^k(t)}
{\sum_{k=0}^{n} p'_k(u_0(t)){u'_0(t)}^k}.$$
Hence the fundamental matrix associated to our variational system is
$$\Phi(t)=\left(\begin{array}{cc}\frac{\sum_{k=0}^{n} 
p_k(u_0)(t){u_0'}^k(t)}{\psi(a)}
&\psi(a)(\sum_{k=0}^{n} p_k(u_0)(t){u_0'}(t)^k)\pi(t)
     \\-\frac{{u_0}'(t)}{\psi(a)}&-\psi(a)u'_0(t)\pi(t)-\psi(a)\pi'(t)
     \frac{\sum_{k=0}^{n} p_k(u_0)(t){u_0'(t)}^k}
{\sum_{k=0}^{n} p'_k(u_0(t)){u_0'(t)}^k}\end{array}\right) .$$
We deduce the principal matrix (the fundamental one with 
$t=\tau_0$).
  $$\Phi(\tau_0)= \left(\begin{array}{cc}1&{\psi(a)}^2 \pi(\tau_0)\\
0&\rho_2\end{array}\right).$$
By the Liouville's formula, we have the characteristic multipliers
 $\rho_1=1$ and $$\rho_2=det(\Phi(\tau_0))=exp(\int_{0}^{\tau_0} (Tr 
{f_x}'(q(\tau))d\tau)=
  exp-(\int_{0}^{\tau_0}\sum_{k=1}^{n} k 
p_k(u_0(\tau)){u'_0(\tau)}^{k-1} )d\tau)$$
  Then we define the  equivalence (\ref{l7}) :
 \begin{equation}\label{la}\rho_2 <1\Longleftrightarrow
  \int_{0}^{\tau_0}(\sum_{k=1}^{n} k p_k(u_0(\tau)){u'_0(\tau)}^{k-1} )d\tau>0
\end{equation}
and the associated Jacobian matrix is : 
$$J(\tau_0)=\left(\begin{array}{cc}-\psi(a)&{\psi(a)}^2 \pi(\tau_0)\\
0&\rho_2-1\end{array}\right).$$
\rightline{$\Box$} 
\subsection{Uniqueness of the periodic solution for an unperturbed equation }
Let us consider now equation (\ref{l4})
$$u''+[u^2+(u+u')^2-1]u'+u=0.$$
That is a special case of generalized Li\'enard equation with $$\varphi(u,u')=(u^2+(u'+u)^2-1)\ and \ \psi(u)=u.$$
We will prove existence and uniqueness of non trivial periodic solution
 for equation (\ref{l4}). Existence will be insured by $C_{LS}$ conditions and for proving
 uniqueness we use a De Castro's result [5] (see also [2]).\\\\
\begin{prop} (De Castro)
Suppose the following system has at least one periodic orbit 
$$\left\{\begin{array} 
{rl}y'&=-\varphi(x,y)y-\psi(x)\\x'&=y.
                                  \end{array}\right.$$ 

Then under the two assumptions, \\
 $a)\;\psi(x)=x;$\\
 $b)\;\varphi(x,y)$ increases, when $|x|$ or $|y|$ or the both increase\\
this periodic orbit is unique.\\
 
\end{prop}
{\text{Let us verify that Equation (6) satisfies the above assumptions}:}\\
\begin{equation}\label{lb}
(\ref{l4})\Longleftrightarrow\left\{\begin{array} {rl} u''+\sum_{k=0}^{3} p_k(u){u'}^k&=0,\\
p_0(u)=\psi(u)=u, \qquad p_1(u)=2u^2-1&, \qquad p_2(u)=2u 
\;\;and \;\;p_3(u)=1.\end{array}
\right.
\end{equation}

Also
 \begin{equation}\label{lc}
(\ref{l4})\Longleftrightarrow \left\{\begin{array} {rl} u''+\varphi(u,u')u'+\psi(u)&=0,\\
\varphi(u,u')=(u^2+(u'+u)^2-1)&\;\;and\;\; \psi(u)=u.\end{array}\right.
    \end{equation}
 
Clearly, the assumptions of Proposition 2 are satified.
In the following, we firstly verify conditions $C_{LS}$ conditions.\\
In that case the 
equation
 $$u''+\varphi(u,u')u'+\psi(u)=0$$ has at least a non trivial periodic 
solution.\\

It is easy to see that $\psi(u)=u$ satisfies $$x \psi(x)>0\qquad 
for \qquad|x|>0$$
$$\int_0^{x} \psi(s) ds =\Psi(x) \quad and \quad lim_{x\rightarrow 
+\infty}
 \Psi(x) = + \infty $$
Now we have $\varphi(0,0)=-1<0.$

By taking $x_0=1, M=1  $ we have $$ \varphi(x,y)\geq 0 
\qquad\;\;\;
 \;\;for\qquad |x|\geq x_0$$ $$\varphi(x,y)\geq -M \qquad for 
\qquad |x|\leq  x_0\\ $$
and the following calculation gives us the optimal value of $x_1>x_0 .$
Let 
\begin{align*}
 H&=\int_{x_0}^{x_1} \varphi(x,y) dx = \int_{1}^{x_1}[ 
x^2+(x+y)^2-1]dx\\
&=\int_{1}^{x_1}[2x^2+2xy+y^2-1]dx = 
[\frac{2}{3}x^3+x^2y+x(y^2-1)]_{1}^{x_1}\\
&=(x_1-1)(\dfrac{{x_1}^2-2x_1+1}{6} 
+2(\frac{x_1+1}{2})^2+2y(\frac{x_1+1}{2})+(y^2-1))\\
&=(x_1-1)(\dfrac{{x_1}^2-2x_1+1}{6} +\varphi(\frac{x_1+1}{2},y))
\end{align*}
Since \ $\frac{x_1+1}{2}\geqslant x_0=1$ \ and using the inequality
$$\varphi(x,y)\geq 0 \qquad\;\;\; \;\;for\qquad |x|\geq x_0 $$
 we then obtain $H\geqslant\frac{(x_1-1)^3}{6}.$\\
Hence, if $\frac{(x_1-1)^3}{6}=10Mx_0=10$, then 
$x_1=1+(60)^\frac{1}{3}$ which satisfies :

$$x_1 > x_0 ,\qquad \int_{x_0}^{x_1} \varphi(x,y) \quad dx\quad \geq 
\quad 10 M x_0 ,$$
for every decreasing function $y(x)>0.$\\
\rightline{$\Box$}
\newpage

\subsection{Existence of periodic solution for perturbed equation 
satisfying $C_{LS}$}
In the following we are dealing with the existence of periodic solution for the equation (\ref{l3}).\\
We assume the initial values :
$$ t=0 ,\;\;u_0(0)=a <0 \;\;and \;\;{u_0}'(0)=0.$$

\begin{theo}\label{montro}
Suppose $1$ is a simple characteristic multiplier of the variational 
system associated to (\ref{l4}).
 Then there are two functions
$\tau,h :U\longrightarrow R$ and constants $\tau_1<\tau_0 $
 such that the periodic solution 
$\nu(t,\alpha,a+h(\epsilon,\alpha),\epsilon,\tau(\epsilon,\alpha))$
  of the equation
 $$ 
u''+{u'}^3+2u{u'}^2+(2u^2-1)u'+u=\epsilon\omega(\frac{t}{\tau},u,u'),$$
 exists for $(\epsilon,\alpha)\in U $ with $|\tau-\tau_0|<\tau_1,$\; 
$\tau(0,0)=\tau_0$ and $h(0,0)=0.$
 \end{theo}

{\textbf{Proof of Theorem \ref{montro}:}}\\
We proceed similarly as in the proof of Theorem 2. We substitute the fundamental matrix, 
 the second characteristic multiplier is $ \rho_2$. 
The following holds for equation (\ref{l4})
$$ \rho_2<1\Longleftrightarrow
  \int_{0}^{\tau_0}(\sum_{k=1}^{3} k p_k(u_0(\tau)){u'_0(\tau)}^{k-1} 
)d\tau>0,$$
 then  $$ \rho_2 <1\Longleftrightarrow 
 \int_{0}^{\tau_0}[2{u_0}^2(\tau)+4u_0(\tau)u'_0(\tau)+3{u'_0(\tau)}^{2}-1]d\tau>0.$$
 It insures that  $1$ is a simple characteristic multiplier of the variational 
system associated to (\ref{l4}) it implies $J(\tau_0)\neq 0.$ Then a periodic solution for 
the perturbed equation (\ref{l3}) exists.

 \rightline{$\Box$}

   \newpage
   Using Scilab we will describe the phase plane of equation (\ref{l4}) $u''+[u^2+(u+u')^2-1]u'+u=0.$
We take
 $x_0=u_0(0)=a=-0.7548829\;,\;y_0={u_0}'(0)=0$\;and the step time of integration $(step=.0001).$\\
 Recall that the periodic orbit is unique.\\

\includegraphics[height=7cm,width=8.25cm]{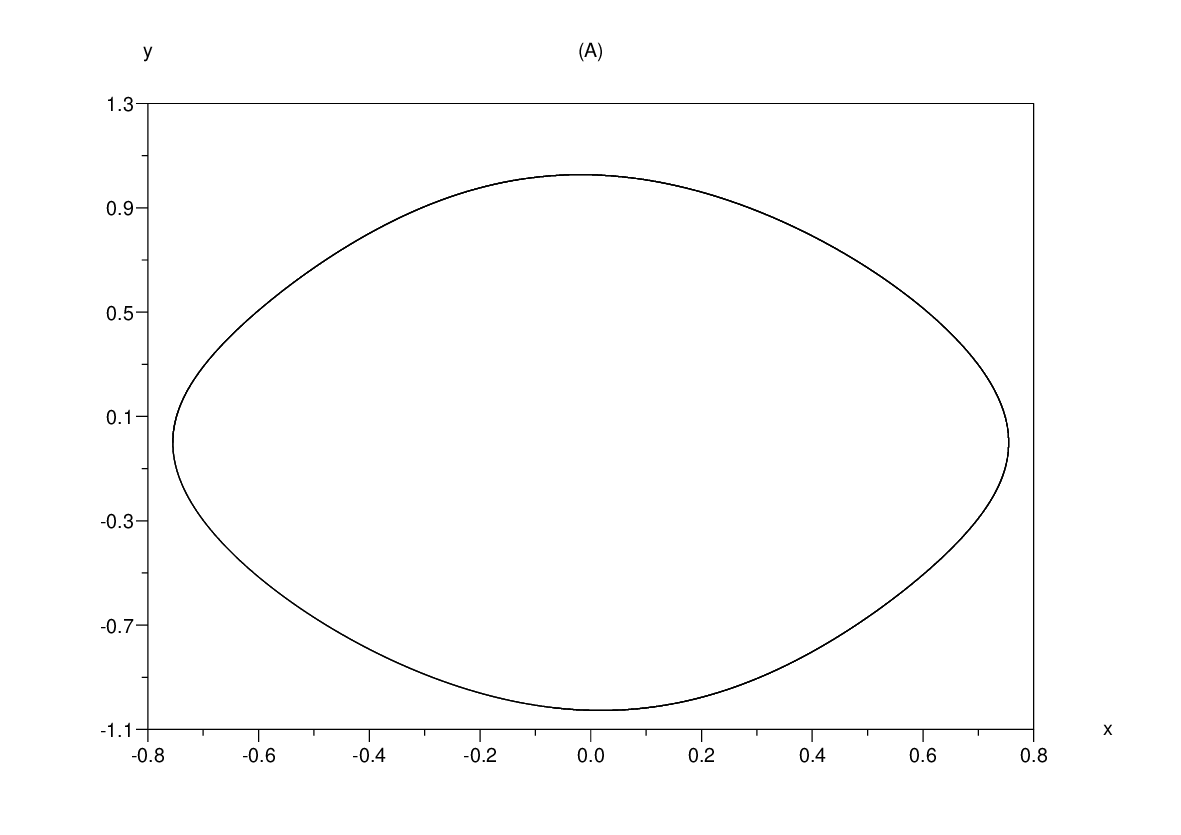}
\includegraphics[height=7cm,width=8.25cm]{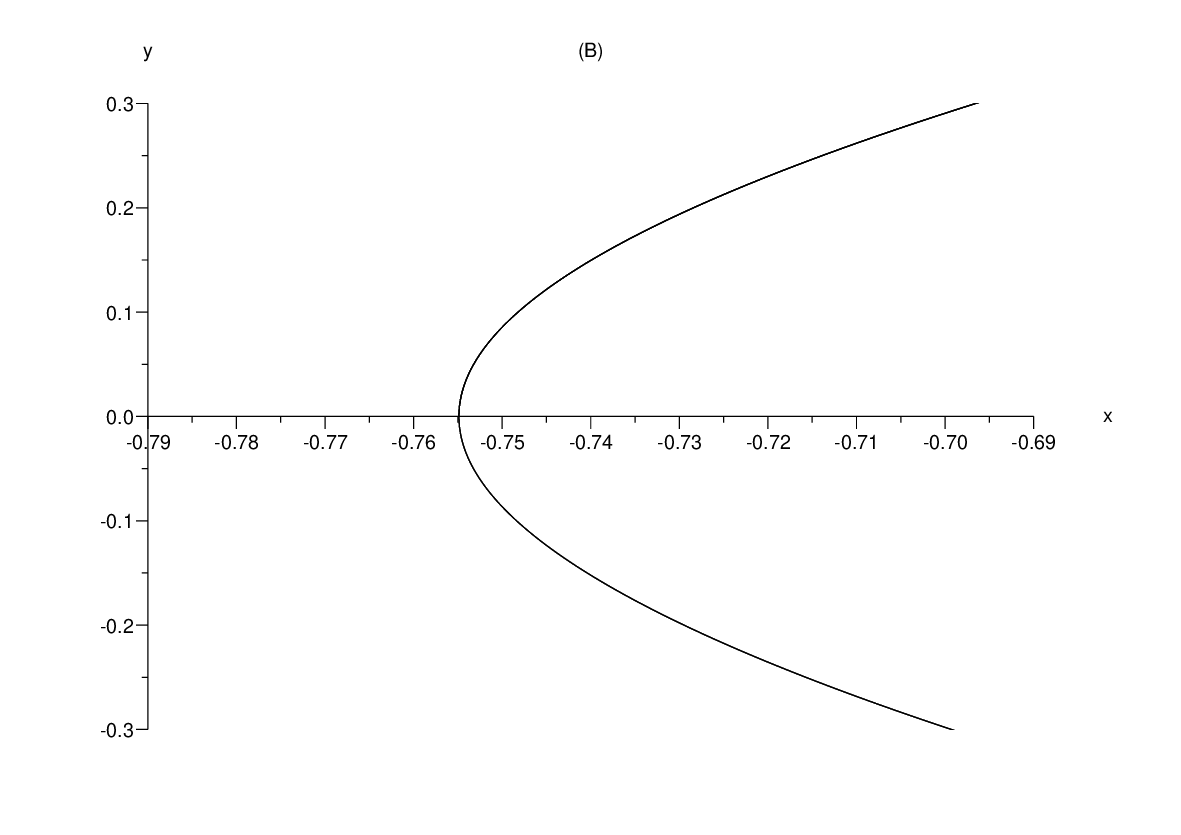}
\qquad\qquad Fig(A)- The unique periodic orbit for the equation $u''+[u^2+(u+u')^2-1]u'+u=0.$
\qquad\qquad Fig(B)- Zoom on the periodic orbit ($\times20$).\\

We take $\epsilon\omega(\frac{t}{\tau},u,u')=\epsilon sin(2t)u'$.
 Some illustrations of the phase portrait for the perturbed equation (\ref{l3}), 
 those can explain existence of a bound $\epsilon_0$, from which periodicity of the
 orbit will be not insured.
 In order to localise $\epsilon_0$, we have taken several values of $\epsilon.$ \\


\includegraphics[height=7cm,width=8.25cm]{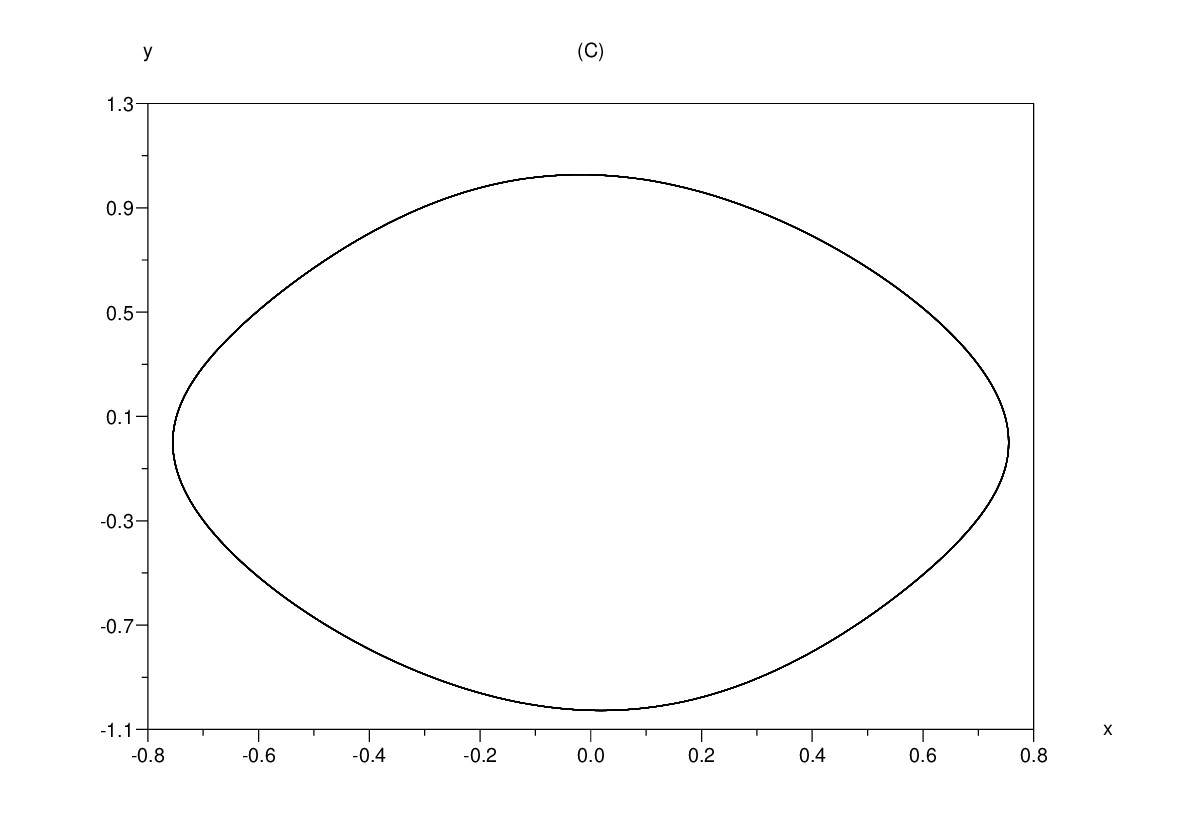}
\includegraphics[height=7cm,width=8.25cm]{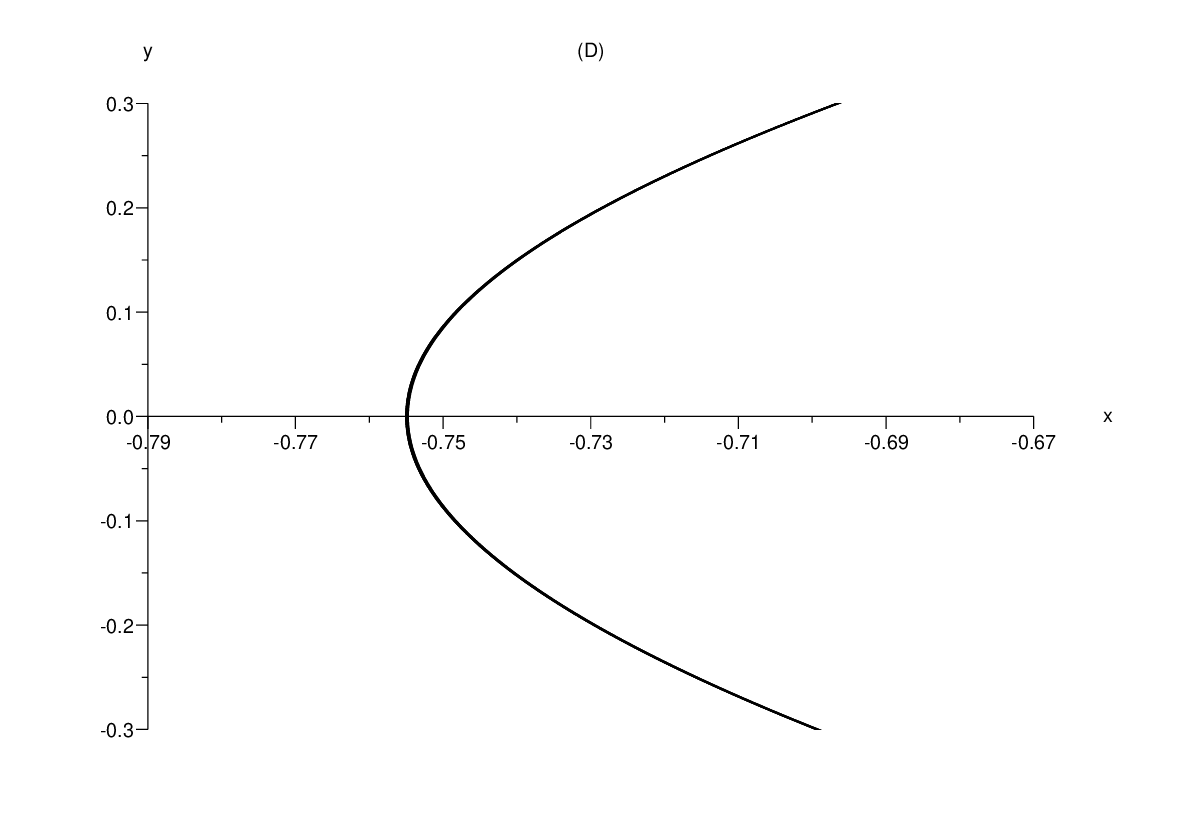}
\qquad\qquad Fig(C)- The periodic orbit for the equation $u''+[u^2+(u+u')^2-1]u'+u=\epsilon\omega(\frac{t}{\tau},u,u')$,
 $\epsilon=0.001$.\\
\qquad\qquad Fig(D)- Zoom on the periodic orbit ($\times20$).

\includegraphics[height=7cm,width=8.25cm]{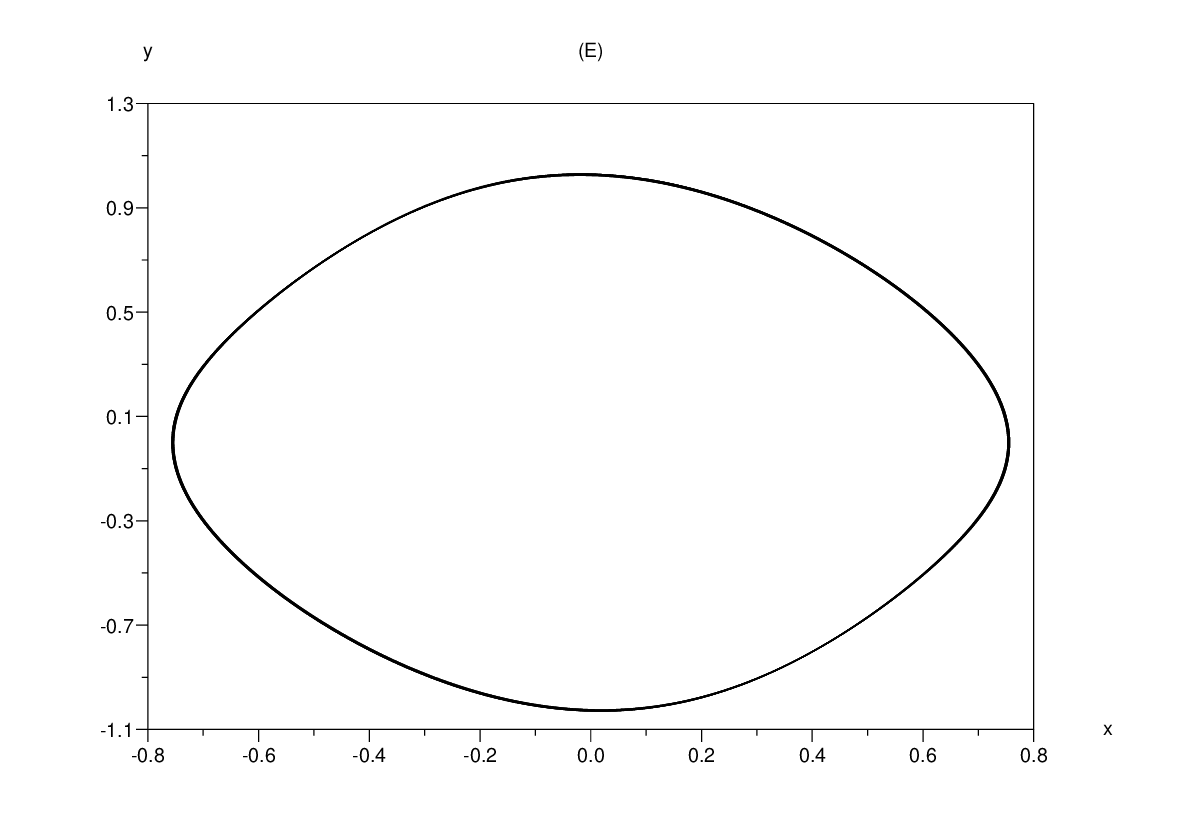}
\includegraphics[height=7cm,width=8.25cm]{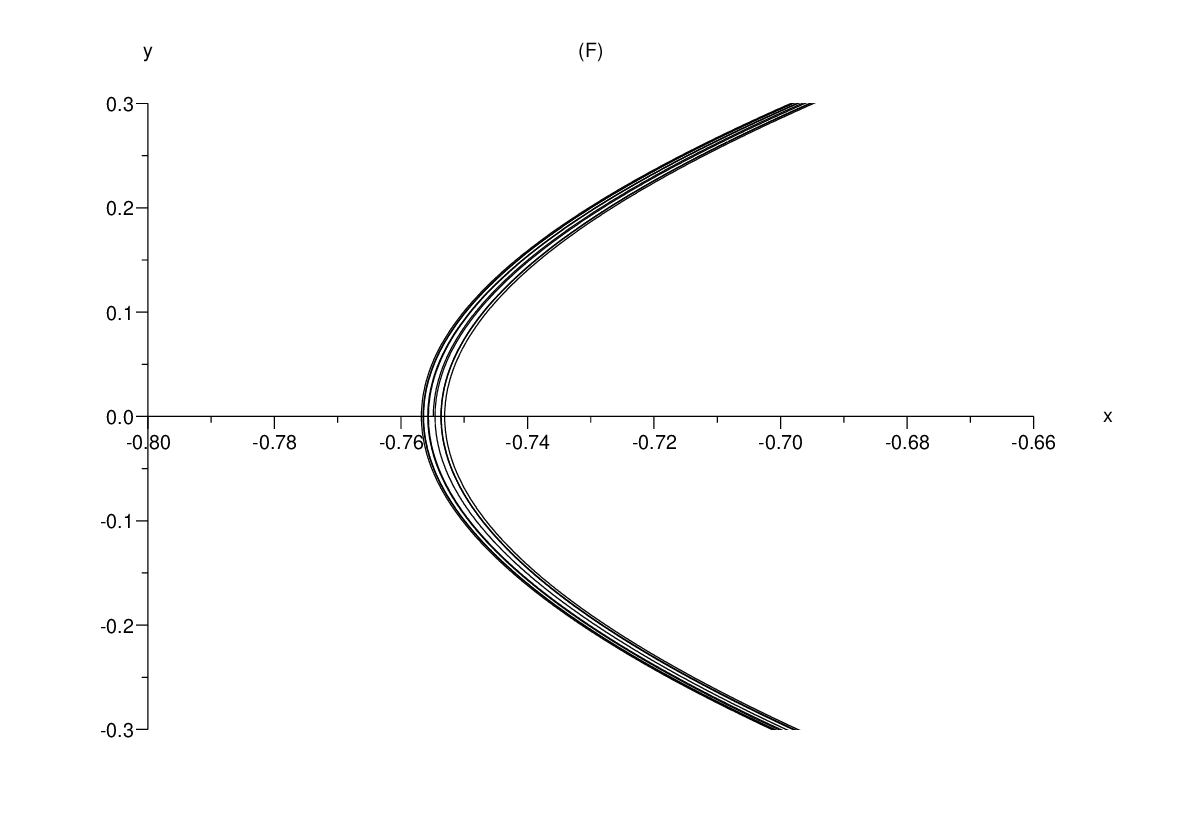}
\qquad\qquad Fig(E)- The orbit for the equation $u''+[u^2+(u+u')^2-1]u'+u=\epsilon\omega(\frac{t}{\tau},u,u')$, $\epsilon=0.01$.\\
\qquad\qquad Fig(F)- Zoom on the orbit ($\times10$) and loss of periodicity.\\
 We see that from the range of $\epsilon=0.01$ the orbit loses the periodicity.\\
For some values of $\epsilon$, we have computed the period .\\

\begin{tabular}{|c|c|c|c|c|c|c|c|c|c|c|}
\hline
$\epsilon$ & 0 & 1/1000 & 1/900 & 1/800 & 1/700 & 1/600 & 1/500 & 1/400 & 1/300 & 1/200 \\
\hline
$\tau$ & 5.4296 & 5.4287 & 5.4286 &5.4285 &5.4283 & 5.4281 & 5.4278 & 5.4274 & 5.4267 & 5.4252 \\
\hline
\end{tabular} \\\\\\\\\\\\\
{\textbf{Acknowledgements}}\\

We thank Professors Miklos Farkas and Jean Marie Strelcyn for their helpful discussions, 
we thank also the referee for the suggestions.



\newpage
  \begin{thebibliography}{20}
\bibitem{C}
        A. R. Chouikha,
        Periodic perturbation of non-conservative second order 
differential equations,
        {\it{Electron.J.Qual.Theory.Differ.Equ}}, 49 (2002), 122-136.
\bibitem{De}
        A.De Castro, 
        Sull'esistenza ed unicità delle soluzioni periodiche dell'equazione\\
         $\ddot x+f(x,\dot x)\dot x+g(x)=0$,
        {\it { Boll. Un. Mat. Ital}}, (3)  9  (1954). 369--372.
\bibitem{F1}
        M. Farkas,
        Periodic motions,
        {\it {Springer-Verlag}},(1994).
\bibitem{L-S}
        N. Levinson and O. K. Smith,
        General equation for relaxation oscillations,\\ {\it{ Duke. Math. Journal.}}
 No 9 (1942), 382-403.
\bibitem{R-S-C}
        R. Reissig  G. Sansonne  R. Conti,
        Qualitative theorie nichtlinearer differentialgleichungen,
	 {\it{Publicazioni del\'l instituto di alta matematica}}, (1963). 
\end{thebibliography}
\end{document}